\documentclass{amsproc}
\usepackage{amssymb}
\usepackage{graphicx}
\usepackage{amscd}
\usepackage{amsmath}
\theoremstyle{plain}

\newtheorem{corollary}{Corollary}

\newtheorem{definition}{Definition}

\newtheorem{lemma}{Lemma}

\newtheorem{proposition}{Proposition}
\newtheorem{remark}{Remark}

\newtheorem{theorem}{Theorem}

\numberwithin{equation}{section}

\begin{document}

\title{Dynamical inverse problem for the discrete Schr\"odinger operator. }
\author{ A. S. Mikhaylov}
\address{St. Petersburg   Department   of   V.A. Steklov    Institute   of   Mathematics
of   the   Russian   Academy   of   Sciences, 7, Fontanka, 191023
St. Petersburg, Russia and Saint Petersburg State University,
St.Petersburg State University, 7/9 Universitetskaya nab., St.
Petersburg, 199034 Russia.} \email{a.mikhaylov@spbu.ru}

\author{ V. S. Mikhaylov}
\address{St.Petersburg   Department   of   V.A.Steklov    Institute   of   Mathematics
of   the   Russian   Academy   of   Sciences, 7, Fontanka, 191023
St. Petersburg, Russia and Saint Petersburg State University,
St.Petersburg State University, 7/9 Universitetskaya nab., St.
Petersburg, 199034 Russia.} \email{v.mikhaylov@spbu.ru}

\date{July, 2016}

\maketitle

{\bf Dedicated to the memory of B. S. Pavlov.}
\newline

\noindent {\bf Abstract.} We consider the inverse problem for the
dynamical system with discrete Schr\"odinger operator and discrete
time. As an inverse data we take a \emph{response operator}, the
natural analog of the dynamical Dirichlet-to-Neumann map. We
derive two types of equations of inverse problem and answer a
question on the characterization of the inverse data, i.e. we
describe the set of operators, which are \emph{response operators}
of the dynamical system governed by the discrete Schr\"odinger
operator.

\section{Introduction.}

The theory of dynamical inverse problems is a wide area of modern
mathematics, by now for all or almost all linear nonstationary
equations of mathematical physics there exist an inverse theory
more or less developed. Theories mostly covers the case of
continuous problems, at the same time just a few attention is paid
to discrete ones. The primary goal of the paper is to improve this
situation.

Let $\mathbb{N}$ be the set of positive natural numbers,
$\mathbb{N}_0=\mathbb{N}\cup\{0\}$. We fix the infinite sequence
of real numbers $(b_1,b_2,\ldots),$ which we call the potential
and consider the dynamical system with discrete time which is a
natural analog of dynamical systems governed by the wave equation
with potential on a semi-axis:
\begin{equation}
\label{Jacobi_dyn} \left\{
\begin{array}l
u_{n,t+1}+u_{n,t-1}-u_{n+1,t}-u_{n-1,t}+b_nu_{n,t}=0,\quad n,t\in \mathbb{N}_0,\\
u_{n,-1}=u_{n,0}=0,\quad n\in \mathbb{N}, \\
u_{0,t}=f_t,\quad t\in \mathbb{N}_0.
\end{array}\right.
\end{equation}
By analogy with continuous problems \cite{B07}, we treat the real
sequence $f=(f_0,f_1,\ldots)$ as a \emph{boundary control}. The
solution to (\ref{Jacobi_dyn}) we denote by $u^f_{n,t}$.

Having fixed $\tau\in \mathbb{N}$, with (\ref{Jacobi_dyn}) we
associate the \emph{response operators}, which maps the control
$f=(f_0,\ldots f_{\tau-1})$ to $u^f_{1,t}$:
\begin{equation}
\left(R^\tau f\right)_t:=u^f_{1,t},\quad t=1,\ldots, \tau.
\end{equation}
The inverse problem we will be dealing with is to recover from
$R^{\tau}$ (part of the) potential $(b_1,b_2,\ldots,b_n)$ for some
$n$. This problems is a natural discrete analog of the inverse
problem for the wave equation where the inverse data is the
dynamical Dirichlet-to-Neumann map, see \cite{B07}.

We will be using the Boundary Control method \cite{B07} which was
initially developed to treat multidimensional dynamical inverse
problems, but since then was applied to multy- and one-
dimensional inverse dynamical, spectral and scattering problems,
problems of signal processing and identification problems
\cite{B13,AMM}.

In the second section we study the forward problem: for
(\ref{Jacobi_dyn}) we prove the analog of d'Alembert integral
representation formula. Prescribing the Dirichlet condition at
$n=N+1$, we consider the second dynamical system with boundary
control at $n=0$ (which will be an analog of the problem on the
finite interval) and develop the solution of this system in
Fourier series. We analyze the dependence of two solutions on the
potential, which lead us to the natural set up of the inverse
problem. In the third section we introduce and prove the
representation formulaes for the main operators of the BC method:
response operator, control and connecting operators. In the forth
section we derive two types of equations for the inverse problem
and give a characterization of the inverse data. In the last
section we point out the connections between the different types
of the inverse data.

The case of the Jacobi matrices of general type as well as the
studying of the inverse spectral problem, i.e. recovering the
semi-infinite matrix from the spectral measure, will be the
subject of forthcoming publications.

\section{Forward problems.}

We fix some positive integer $T$. By $\mathcal{F}^T$ we denote the
space of controls: $\mathcal{F}^T:=\mathbb{R}^T$, $f\in
\mathcal{F}^T$, $f=(f_0,\ldots,f_{T-1})$.

First, we derive the representation formulas for the solution to
(\ref{Jacobi_dyn}) which could be considered as analogs of known
formulas for the wave equation \cite{AM}.
\begin{lemma}
The solution to (\ref{Jacobi_dyn}) admits the representation
\begin{equation}
\label{Jac_sol_rep}
u_{n,t}=f_{t-n}+\sum_{s=n}^{t-1}w_{n,s}f_{t-s-1},\quad n,t\in
\mathbb{N}_0.
\end{equation}
where $w_{n,s}$ satisfies the Goursat problem
\begin{equation}
\label{Goursat} \left\{
\begin{array}l
w_{n,t+1}+w_{n,t-1}-w_{n+1,t}-w_{n-1,t}+b_nw_{n,t}=0,\quad n,s\in \mathbb{N}_0, \,\,s>n,\\
w_{n,n}=-\sum_{k=1}^n b_k,\quad n\in \mathbb{N},\\
w_{0,t}=0,\quad t\in \mathbb{N}_0.
\end{array}
\right.
\end{equation}
\end{lemma}
\begin{proof}
We assume that $u^f_{n,t}$ has a form (\ref{Jac_sol_rep}) with
unknown $w_{n,s}$ and plug it to equation in (\ref{Jacobi_dyn}):
\begin{eqnarray*}
0=b_nf_{t-n}+\sum_{s=n}^{t-1}b_nw_{n,s}f_{t-s-1}+\sum_{s=n}^{t}w_{n,s}f_{t-s}+\sum_{s=n}^{t-2}w_{n,s}f_{t-s-2}\\
-\sum_{s=n+1}^{t-1}w_{n,s}f_{t-s-1}-\sum_{s=n-1}^{t-1}w_{n-1,s}f_{t-s-1}.
\end{eqnarray*}
Changing the order of summation we get
\begin{eqnarray*}
0=b_nf_{t-n}+w_{n+1,n}f_{t-n-1}
-w_{n-1,n-1}f_{t-n}+\sum_{s=n}^{t-1}f_{t-s-1}\left(b_nw_{n,s}-w_{n+1,s}-w_{n-1,s}\right)\\
+\sum_{s=n-1}^{t-1}w_{n,s+1}f_{t-s-1}+\sum_{s=n+1}^{t-1}w_{n,s-1}f_{t-s-1}=f_{t-n-1}(w_{n+1,n}-w_{n,n-1})+b_nf_{t-n}\\
+\sum_{s=n}^{t-1}f_{t-s-1}\left(w_{n,s+1}+w_{n.s-1}-w_{n+1,s}-w_{n-1,s}+b_nw_{n,s}\right)+f_{t-n}(w_{n,n}-w_{n-1,n-1}).
\end{eqnarray*}
Counting that $w_{n,s}=0$ when $n>s$ and arbitrariness of $f\in
\mathcal{F}^T$, we arrive at (\ref{Jac_sol_rep}).
\end{proof}

We fix $N\in \mathbb{N}$. Along with (\ref{Jacobi_dyn}) we
consider the analog of the wave equation with the potential on the
interval: we assume that $(b_n)$ is finite: $n=1,\ldots,N$ and
impose the Dirichlet condition at $n=N+1$. Then for a control
$f=(f_0,f_1,\ldots)$ we consider
\begin{equation}
\label{Jacobi_dyn_int} \left\{
\begin{array}l
v_{n,t+1}+v_{n,t-1}-v_{n+1,t}-v_{n-1,t}+b_nv_{n,t}=0,\quad t\in \mathbb{N}_0,\,\, n\in 0,\ldots, N+1\\
v_{n,-1}=v_{n,0}=0,\quad n=1,2,\ldots,N+1 \\
v_{0,t}=f_t,\quad v_{N+1,t}=0,\quad t\in \mathbb{N}_0.
\end{array}\right.
\end{equation}
We denote the solution to (\ref{Jacobi_dyn_int}) by $v^f$.

Let $\phi_n(\lambda)$ be the solution to
\begin{equation}
\label{Spec_sol} \left\{
\begin{array}l \phi_{i+1}+\phi_{i-1}-b_n\phi_i=\lambda\phi_i,\\
\phi_0=0,\,\,\phi_1=1.
\end{array}
\right.
\end{equation}
We introduce the Hamiltonian
\begin{equation*}
H_N:=\begin{pmatrix} -b_1 & 1 & 0 & \ldots & 0\\
1 & -b_2 & 1 & \ldots & 0\\
\cdot & \cdot &\cdot &\cdot &\cdot \\
0 &\ldots & 0& 1& -b_N
\end{pmatrix}
\end{equation*}
Let $\{\varphi^k,\lambda_k\}_{k=1}^N$ be eigenvectors chosen such
that $\varphi^k_1=1$ and eigenvalues of $H_N$. Introduce the
numbers $\rho_k$ by
\begin{equation}
\label{Ortog} (\varphi^k,\varphi^l)=\delta_{kl}\rho_k,
\end{equation}
where $(\cdot,\cdot)$-- is a scalar product in $\mathbb{R}^N$.
\begin{definition}
The set
\begin{equation}
\label{SP_data} \{\lambda_k,\rho_k\}_{k=1}^N
\end{equation}
is called the spectral data.
\end{definition}
On introducing vectors $\phi^n\in \mathbb{R}^N$ by the rule
$\phi^n_i:=\phi_i(\lambda_n)$, $n,i=1,\ldots,N,$ we have
\begin{proposition}
The solutions of $\phi_{N+1}(\lambda)=0$ are $\lambda_n$,
$n=1,\ldots,N$; and $\phi^n_i=\varphi^n_i$, $n,i=1,\ldots,N.$
\end{proposition}
\begin{proof}
Take eigenvector $\varphi^n$ corresponding to eigenvalue
$\lambda_n$ and compare it with $\phi^n$. By the definition of
$\varphi^n$ and condition in (\ref{Spec_sol}):
$\varphi^n_1=\phi^n_1=1$. On the other hand, comparing the first
line in the equation on eigenvalues
$H_N\varphi^n=\lambda^n\varphi^n$ and (\ref{Spec_sol}) for $i=1$
we have:
\begin{eqnarray*}
-\varphi^n_1b_1+\varphi^n_2=\lambda_n\varphi^n_1,\\
\phi^n_2-b_1\phi^n_1=\lambda_n\phi^n_1,
\end{eqnarray*}
which implies $\varphi^n_2=\phi^n_2$, for $k<N$ comparing $k-$th
line in $H_N\varphi^n=\lambda^n\varphi^n$ and (\ref{Spec_sol}) for
$i=k$, we arrive at $\varphi^n_{k+1}=\phi^n_{k+1}$. And for $k=N:$
\begin{eqnarray*}
-\varphi^n_Nb_N+\varphi^n_{N-1}=\lambda_n\varphi^n_N,\\
\phi^n_{N+1}+\phi^n_{N-1}-b_N\phi^n_N=\lambda_n\phi^n_N,
\end{eqnarray*}
which holds if and only if $\phi^n_{N+1}(\lambda_n)=0.$
\end{proof}

We take $y\in \mathbb{R}^N, \,\,y=(y_1,\ldots,y_N)$, for each $n$
we multiply the equation in (\ref{Jacobi_dyn_int}) by $y_n$, sum
up and evaluate the following expression, changing the order of
summation
\begin{eqnarray}
0=\sum_{n=1}^N \left(v_{n,t+1}y_n+v_{n,t-1}y_n - v_{n+1,t}y_n-
v_{n-1,t}y_n+b_nv_{n,t}y_n\right)\notag\\
=\sum_{n=1}^N \left(v_{n,t+1}y_n+v_{n,t-1}y_n -
v_{n,t}(y_{n-1}+y_{n+1})+b_nv_{n,t}y_n\right)\notag\\
-v_{N+1,t}y_N-v_{0,t}y_1+v_{1,t}y_0+v_{N,t}y_{N+1}\label{int_rav}
\end{eqnarray}
Now we choose $y=\varphi^l$, $l=1\ldots, N$. On counting that
$\varphi^l_0=\varphi^l_{N+1}=0,$ $\varphi^l_1=1,$ $v_{0,t}=f_t$,
$v_{N+1,t}=0$ we evaluate (\ref{int_rav}) arriving at:
\begin{equation}
\label{Integr_eqn} 0=\sum_{n=1}^N
\left(v_{n,t+1}\varphi^l_n+v_{n,t-1}\varphi^l_n -
v_{n,t}\left(\varphi^l_{n-1}+\varphi^l_{n+1}-b_n\varphi^l_n\right)\right)
-f_t=0
\end{equation}

\begin{definition}
For $a,b\in l^\infty$ we define the convolution $c=a*b\in
l^\infty$ by the formula
\begin{equation*}
c_t=\sum_{s=0}^{t}a_sb_{t-s},\quad t\in \mathbb{N}
\end{equation*}
\end{definition}

We assume that the solution to (\ref{Jacobi_dyn_int}) has a form
\begin{equation}
\label{Jac_sol_rep_int} v^f_{n,t}= \left\{
\begin{array}l
\sum_{k=1}^N c_t^k\varphi^k_n,\quad n=1,\ldots,N\\
f_t,\quad n=0.
\end{array}
\right.
\end{equation}
\begin{proposition}
The coefficients $c^k$ admits the representation:
\begin{equation}
\label{c_koeff_repr}
c^k=\frac{1}{\rho_k}T\left(\lambda_k\right)*f,
\end{equation}
where
$T(2\lambda)=(T_1(2\lambda),T_2(2\lambda),T_3(2\lambda),\ldots)$
are Chebyshev polynomials of the second kind.
\end{proposition}
\begin{proof}
We plug (\ref{Jac_sol_rep_int}) into (\ref{Integr_eqn}) and
evaluate, counting that
$\varphi^l_{n-1}+\varphi^l_{n+1}-b_n\varphi^l_n=\lambda_l\varphi^l_n$:
\begin{eqnarray*}
\sum_{n=1}^N \left(v_{n,t+1}+v_{n,t-1}-\lambda_l
v_{n,t}\right)\varphi^l_n =f_t,\\
\sum_{n=1}^N
\sum_{k=1}^N\left(c^k_{t+1}\varphi^k_{n}+c^k_{t-1}\varphi^k_{n}-\lambda_l
c^k_{t}\varphi^k_{n}\right)\varphi^l_n =f_t.
\end{eqnarray*}
Changing the order of summation and using (\ref{Ortog}) we finally
arrive at the following  equation on $c^k_t$, $k=1,\ldots, N$:
\begin{equation}
\label{c_eqn} \left\{
\begin{array}l
c^k_{t+1}+c^k_{t-1}-\lambda_k c^k_{t}=\frac{1}{\rho_k}f_t,\\
c^k_{-1}=c^k_0=0.
\end{array}
\right.
\end{equation}
We assume that solution to (\ref{c_eqn}) has a form
$c^k=\frac{1}{\rho_k}T*f,$ or
\begin{equation}
\label{c_eqn_1} c^k_t=\frac{1}{\rho_k}\sum_{l=0}^{t}T_lf_{t-l}.
\end{equation}
Plugging it into (\ref{c_eqn}), we get
\begin{eqnarray*}
\frac{1}{\rho_k}\left(\sum_{l=0}^{t+1}f_lT_{t+1-l}+
\sum_{l=0}^{t-1}f_l T_{t-1-l}-\lambda_k\sum_{l=0}^{t}f_l T_{t-l}\right)=\frac{1}{\rho_k}f_t\\
\sum_{l=0}^{t}f_l\left(T_{t+1-l}+T_{t-1-l}-\lambda_k
T_{t-l}\right)+f_t T_1-f_{t-1}T_0 =f_t.
\end{eqnarray*}
We see that (\ref{c_eqn_1}) holds if $T$ solves
\begin{equation*}
\left\{
\begin{array}l
T_{t+1}+T_{t-1}-\lambda_k T_{t}=0,\\
T_{0}=0,\,\, T_1=1.
\end{array}
\right.
\end{equation*}
Thus $T_k(2\lambda)$ are Chebyshev polynomials of the second kind.

\end{proof}

\section{Operators of the the BC method.}


As an inverse data for (\ref{Jacobi_dyn}) we use the analog of the
dynamical response operator (dynamical Dirichlet-to-Neumann map)
\cite{B07}.
\begin{definition}
For (\ref{Jacobi_dyn}) the \emph{response operator}
$R^T:\mathcal{F}^T\mapsto \mathbb{R}^T$ is defined by the rule
\begin{equation*}
\left(R^Tf\right)_t=u^f_{1,t}, \quad t=1,\ldots,T.
\end{equation*}
\end{definition}
Introduce the notation: the \emph{response vector} is the
convolution kernel of the response operator,
$r=(r_0,r_1,\ldots,r_{T-1})=(1,w_{1,1},w_{1,2},\ldots w_{1,T-1})$.
Then in accordance with (\ref{Jac_sol_rep})
\begin{eqnarray}
\label{R_def}
\left(R^Tf\right)_t=u^f_{1,t}=f_{t-1}+\sum_{s=1}^{t-1}
w_{1,s}f_{t-1-s}
\quad t=1,\ldots,T.\\
\notag \left(R^Tf\right)=r*f_{\cdot-1},\quad \text{where $r_0=1$}.
\end{eqnarray}

For the system (\ref{Jacobi_dyn_int}) we introduce the response
operator by
\begin{definition}
For the system in (\ref{Jacobi_dyn_int}) the \emph{response
operator} $R^T_i:\mathcal{F}^T\mapsto \mathbb{R}^T$ is defined by
the rule
\begin{equation}
\label{R_def_int} \left(R^T_if\right)_t=v^f_{1,t}, \quad
t=1,\ldots,T.
\end{equation}
\end{definition}
The corresponding \emph{response vector} we denote by
$(r^i_1,r^i_2,\ldots)$. More information on this operator and on
the inverse spectral problem one can find in the last section.

We introduce the inner space of dynamical system
(\ref{Jacobi_dyn}) $\mathcal{H}^T:=\mathbb{R}^T$, $h\in H^T$,
$h=(h_1,\ldots, h_T)$. For (\ref{Jacobi_dyn}) The \emph{control
operator} $W^T:\mathcal{F}^T\mapsto \mathcal{H}^T$ is defined by
the rule
\begin{equation*}
W^Tf:=u^f_{n,T},\quad n=1,\ldots,T.
\end{equation*}
Directly from (\ref{Jac_sol_rep}) we deduce that
\begin{equation}
\label{W^T_rep}
\left(W^Tf\right)_n=u^f_{n,T}=f_{T-n}+\sum_{s=n}^{T-1}w_{n,s}f_{T-s-1},\quad
n=1,\ldots,T.
\end{equation}
The following statement imply the controllability of the dynamical
system (\ref{Jacobi_dyn}).
\begin{theorem}
\label{teor_control} The operator $W^T$ is an isomorphism between
$\mathcal{F}^T$ and $\mathcal{H}^T$.
\end{theorem}
\begin{proof}
We fix some $a\in \mathcal{H}^T$ and look for a control $f\in
\mathcal{F}^T$ such that $W^Tf=a$. To this aim we write down the
operator as
\begin{equation}
W^Tf=\begin{pmatrix} u_{1,T}\\
u_{2,T}\\
\cdot\\
u_{k,T}\\
\cdot\\
u_{T,T}
\end{pmatrix}=\begin{pmatrix}
1 & w_{1,1} & w_{1,2} & \ldots & \ldots & w_{1,T-1}\\
0 & 1 & w_{2,2} & \ldots & \ldots & w_{2, T-1}\\
\cdot & \cdot & \cdot & \cdot & \cdot & \cdot \\
0 & \ldots & 1& w_{k,k} & \ldots & w_{k, T-1}\\
\cdot & \cdot & \cdot & \cdot & \cdot & \cdot \\
0 & 0 & 0 & 0 & \ldots & 1
\end{pmatrix}
\begin{pmatrix} f_{T-1}\\
f_{T-2}\\
\cdot\\
f_{T-k-1}\\
\cdot\\
f_{0}
\end{pmatrix}
\end{equation}
We introduce the notations
\begin{eqnarray*}
J_T: \mathcal{F}^T\mapsto \mathcal{F}^T,\quad
\left(J_Tf\right)_n=f_{T-1-n},\quad n=0,\ldots, T-1,\\
K\in \mathbb{R}^{T\times T},\quad k_{ij}=0,\, i>j, k_{ii}=1,\,
k_{ij}=w_{ij-1},\,i<j.
\end{eqnarray*}
Then $W^T=\left(I+K\right)J^T$. Obviously, this operator is
invertible, which proves the statement of the theorem.
\end{proof}

For the system (\ref{Jacobi_dyn_int}) the \emph{control operator}
$W^T_i:\mathcal{F}^T\mapsto \mathcal{H}^N$ is defined by the rule
\begin{equation*}
W^T_if:=v^f_{n,T},\quad n=1,\ldots,N.
\end{equation*}
The representation for this operator immediately follows from
(\ref{Jac_sol_rep_int}), (\ref{c_koeff_repr}).

For the system (\ref{Jacobi_dyn}) we introduce the
\emph{connecting operator} $C^T: \mathcal{F}^T\mapsto
\mathcal{F}^T$ by the quadratic form: for arbitrary $f,g\in
\mathcal{F}^T$ we define
\begin{equation}
\label{C_T_def} \left(C^T
f,g\right)_{\mathcal{F}^T}=\left(u^f_{\cdot,T},
u^g_{\cdot,T}\right)_{\mathcal{H}^T}=\left(W^Tf,W^Tg\right)_{\mathcal{H}^T}.
\end{equation}
We observe that $C^T=\left(W^T\right)^*W^T$, so $C^T$ is an
isomorphism in $\mathcal{F}^T$. The fact that $C^T$ can be
expressed in terms of response $R^{2T}$ is crucial in BC-method.
\begin{theorem}
Connecting operator admits the representation in terms of inverse
data:
\begin{equation}
\label{C_T_repr} C^T=C^T_{ij},\quad
C^T_{ij}=\sum_{k=0}^{T-\max{i,j}}r_{|i-j|+2k},\quad r_0=1.
\end{equation}
\begin{equation*}
C^T=
\begin{pmatrix}
1+r_2+\ldots+r_{2T-2} & r_1+r_3+\ldots+r_{2T-3} & \ldots &
r_T+r_{T-2} &
r_{T-1}\\
r_1+r_3+\ldots+r_{2T-3} & 1+r_2+\ldots+r_{2T-4} & \ldots & \ldots
&r_{T-2}\\
\cdot & \cdot & \cdot & \cdot & \cdot \\
r_{T-3}+r_{T-1}+r_{T+1} &\ldots & 1+r_2+r_4 & r_1+r_3 & r_2\\
r_{T}+r_{T-2}&\ldots &r_1+r_3&1+r_2&r_1 \\
r_{T-1}& r_{T-2}& \ldots & r_1 &1
\end{pmatrix}
\end{equation*}
\end{theorem}
\begin{proof}
For fixed $f,g\in \mathcal{F}^T$ we introduce the
\emph{Blagoveshchensky function} by
\begin{equation*}
\psi_{n,t}:=\left(u^f_{\cdot,n},
u^g_{\cdot,t}\right)_{\mathcal{H}^T}=\sum_{k=1}^T
u^f_{k,n}u^g_{k,t}
\end{equation*}
Then we show that $\psi_{n,t}$ satisfies some difference equation.
Indeed, we can evaluate:
\begin{eqnarray*}
\psi_{n,t+1}+\psi_{n,t-1}-\psi_{n+1,t}-\psi_{n-1,t}=\sum_{k=1}^T
u^f_{k,n}\left(u^g_{k,t+1}+u^g_{k,t-1}\right)\\
-\sum_{k=1}^T \left(u^f_{k,n+1}+u^f_{k,n-1}\right)u^g_{k,t}=
\sum_{k=1}^T u^f_{k,n}\left(u^g_{k+1,t}+u^g_{k-1,t}\right)
\\
-\sum_{k=1}^T
\left(u^f_{k+1,n}+u^f_{k-1,n}\right)u^g_{k,t}=\sum_{k=1}^T
u^g_{k,t}
\left(u^f_{k+1,n}+u^f_{k-1,n}\right)+u^g_{0,t}u^f_{1,n}-u^f_{0,n}u^g_{1,t}\\
+u^g_{T+1,t}u^f_{T,n}-u^f_{T+1,n}u^g_{T,t}-\sum_{k=1}^T u^g_{k,t}
\left(u^f_{k+1,n}+u^f_{k-1,n}\right)=g_t(Rf)_n-f_n(Rg)_t
\end{eqnarray*}
So we arrive at the following boundary problem for $\psi_{n,t}$:
\begin{eqnarray}
\label{Blag_eqn}
&\left\{
\begin{array}l
\psi_{n,t+1}+\psi_{n,t-1}-\psi_{n+1,t}-\psi_{n-1,t}=h_{n,t},\quad n,t\in \mathbb{N}_0,\\
\psi_{0,t}=0,\,\, \psi_{n,0}=0,
\end{array}
\right. \\
&h_{n,t}=g_t(Rf)_n-f_n(Rg)_t\notag
\end{eqnarray}
We introduce the set
\begin{eqnarray*}
K(n,t):=\left\{(n,t)\cup \{(n-1,t-1),
(n+1,t-1)\}\cup\{(n-2,t-2),(n,t-2),\right.\\
\left.(n+2,t-2)\}\cup\ldots\cup\{(n-t,0),(n-t+2,0),\ldots,(n+t-2,0),(n+t,0)\}\right\}\\
=\bigcup_{\tau=0}^t\bigcup_{k=0}^\tau
\left(n-\tau+2k,t-\tau\right).
\end{eqnarray*}
The solution to (\ref{Blag_eqn}) is given by
\begin{equation*}
\psi_{n,t}=\sum_{k,\tau\in K(n,t-1)}h(k,\tau).
\end{equation*}
We observe that $\psi_{T,T}=\left(C^Tf,g\right)$, so
\begin{equation}
\label{C_T_sol}\left(C^Tf,g\right)=\sum_{k,\tau\in
K(T,T-1)}h(k,\tau).
\end{equation}
Notice that in the r.h.s. of (\ref{C_T_sol}) the argument $k$ runs
from $1$ to $2T-1.$ We extend $f\in \mathcal{F}^T$,
$f=(f_0,\ldots,f_{T-1})$ to $f\in \mathcal{F}^{2T}$ by:
\begin{equation*}
f_T=0,\quad f_{T+k}=-f_{T-k},\,\, k=1,2,\ldots, T-1.
\end{equation*}
Due to this odd extension, $\sum_{k,\tau\in K(T,T-1)}
f_k(R^Tg)_\tau=0$, so (\ref{C_T_sol}) gives
\begin{eqnarray*}
\left(C^Tf,g\right)=\sum_{k,\tau\in K(T,T-1)}
g_\tau\left(R^{2T}f\right)_k=g_0\left[\left(R^{2T}f\right)_1+\left(R^{2T}f\right)_3+\ldots+\left(R^{2T}f\right)_{2T-1}\right]\\
+g_1\left[\left(R^{2T}f\right)_2+\left(R^{2T}f\right)_4+\ldots+\left(R^{2T}f\right)_{2T-2}\right]+\ldots+g_{T-1}\left(R^{2T}f\right)_{T}.
\end{eqnarray*}
Finally we infer that
\begin{equation*}
C^Tf=\left(\left(R^{2T}f\right)_1+\ldots+\left(R^{2T}f\right)_{2T-1},\left(R^{2T}f\right)_2+\ldots+\left(R^{2T}f\right)_{2T-2},\ldots,\left(R^{2T}f\right)_{T}
\right)
\end{equation*}
from where the statement of the theorem follows.
\end{proof}
One can observe that $C^T_{ij}$ satisfies the difference boundary
problem.
\begin{corollary}
The kernel of $C^T$ satisfy
\begin{equation*}
\left\{
\begin{array}l
C^T_{i,j+1}+C^T_{i,j-1}-C^T_{i+1,j}-C^T_{i-1,j}=0,\\
C^T_{i,T}=r_{T-i},\,\,C^T_{T,j}=r_{T-j},\,\, r_0=1.
\end{array}
\right.
\end{equation*}
\end{corollary}

For the system (\ref{Jacobi_dyn_int}) the \emph{connecting
operator} $C^T_i: \mathcal{F}^T\mapsto \mathcal{F}^T$ is
introduced in the similar way: for arbitrary $f,g\in
\mathcal{F}^T$ we define
\begin{equation}
\label{C_T_def_int} \left(C^T_i
f,g\right)_{\mathcal{F}^T}=\left(v^f_{\cdot,T},
v^g_{\cdot,T}\right)_{\mathcal{H}^N}=\left(W^T_if,W^T_ig\right)_{\mathcal{H}^N}.
\end{equation}
More information on $C^T_i$ one can find in the final section.

\section{Inverse problem. }

The dependence of the solution (\ref{Jacobi_dyn}) $u^f$ on the
potential $(b_1,b_2,\ldots)$ resemble one of the wave equation
with the potential: take some $M\in \mathbb{N}$. From the very
equation one can see that the term $u^f_{n,t}$ with smallest
$\{n,t\}$, which depends on $b_M$ is $u^f_{M,M+1}$. Thus,
$u^f_{1,t}$ become to depend on $b_M$ starting from $t=2M$. This
is an analog of the effect of the finite speed of wave propagation
in the wave equation. Consider (\ref{Jacobi_dyn_int}) with $N=M$.
We see that the solution to (\ref{Jacobi_dyn_int}) $v^f_{1,t}$
does not "feel" the boundary condition at $n=M+1$:
$u^f_{1,t}=v^f_{1,t}$ for $t=1,\ldots,2M$. Or in other words that
means that $R^{2M} = R_i^{2M}$. This leads to the following
natural set up of the inverse problem: By the given operator
$R^{2M}$ to recover the (part) of the potential
$(b_1,\ldots,b_{M})$. In what follows we will be dealing with the
IP for the system (\ref{Jacobi_dyn}), only in the last section we
comment on the system (\ref{Jacobi_dyn_int}).

\subsection{Krein equations}
Let $\alpha,\beta\in \mathbb{R}$ and $y$ be solution to
\begin{equation}
\label{y_special} \left\{
\begin{array}l
y_{k+1}+y_{k-1}-b_ky_k=0,\\
y_0=\alpha,\,\, y_1=\beta.
\end{array}
\right.
\end{equation}
We set up the following control problem: to find a control $f^T\in
\mathcal{F}^T$ such that
\begin{equation}
\label{Control_probl}
\left(W^Tf^T\right)_k=y_k,\quad
k=1,\ldots,T.
\end{equation}
Due to Theorem \ref{teor_control}, this problem has unique
solution. Let $\varkappa^T$ be a solution to
\begin{equation}
\label{kappa} \left\{
\begin{array}l
\varkappa^T_{t+1}+\varkappa^T_{t-1}=0,\quad t=0,\ldots,T,\\
\varkappa^T_{T}=0,\,\, \varkappa^T_{T-1}=1.
\end{array}
\right.
\end{equation}
We show that the control $f^T$ satisfies the Krein equation:
\begin{theorem}
The control $f^T$, defined by (\ref{Control_probl}) satisfies the
following equation in $\mathcal{F}^T$:
\begin{equation}
\label{C_T_Krein} C^Tf^T=\beta\varkappa^T-\alpha
\left(R^T\right)^*\varkappa^T.
\end{equation}
\end{theorem}
\begin{proof}
Let us take $f^T$ solving (\ref{Control_probl}). We observe that
for any fixed $g\in \mathcal{F}^T$:
\begin{equation}
\label{Kr_1}
u^g_{k,T}=\sum_{t=1}^{T-1}\left(u^g_{k,t+1}+u^g_{k,t-1}\right)\varkappa^T_t.
\end{equation}
Indeed, changing the order of summation in the r.h.s. of
(\ref{Kr_1}), we get
\begin{equation*}
\sum_{t=1}^{T-1}\left(u^g_{k,t+1}+u^g_{k,t-1}\right)\varkappa^T_t=\sum_{t=1}^{T-1}\left(\varkappa^T_{t+1}+\varkappa^T_{t-1}
\right)u^g_{k,t}+u^g_{k,0}\varkappa^T_1-u^g_{k,T}\varkappa^T_{T-1}.
\end{equation*}
which gives (\ref{Kr_1}) due to (\ref{kappa}). Using this
observation, we can evaluate
\begin{eqnarray*}
\left(C^Tf^T,g\right)=\sum_{k=1}^T
y_ku^g_{k,T}=\sum_{k=1}^T\sum_{t=0}^{T-1}\left(u^g_{k,t+1}+u^g_{k,t-1}\right)\varkappa^T_t
y_k\\
=\sum_{t=0}^{T-1}\varkappa^T_t\left(\sum_{k=1}^T
\left(u^g_{k+1,t}y_k+u^g_{k-1,t}y_k - b_ku^g_{k,t}y_k\right)\right)\\
=\sum_{t=0}^{T-1}\varkappa^T_t\left(\sum_{k=1}^T
\left(u^g_{k,t}(y_{k+1}+y_{k-1}-b_ky_k\right)+u^g_{0,t}y_1+u^g_{T+1,t}y_T-u^g_{1,t}y_0-u^g_{T,t}y_{T+1}
\right)\\
=\sum_{t=0}^{T-1}\varkappa^T_t\left(\beta
g_t-\alpha\left(R^Tg\right)_t \right)=\left(\varkappa^T, \beta
g-\alpha \left(R^Tg\right)\right)=\left(\beta\varkappa^T - \alpha
\left(\left(R^T\right)^*\varkappa^T\right), g\right).
\end{eqnarray*}
From where (\ref{C_T_Krein}) follows.
\end{proof}

Having found $f^\tau$ for $\tau=1,\ldots,T$, we can recover the
potential $b_n,$ $n=1,\ldots,T-1$. Indeed: by the constructions of
$f^\tau$ we have $\left(W^\tau f^\tau\right)_{\tau}=y_\tau$, on
the other hand, from (\ref{W^T_rep}) we can infer that
$\left(W^\tau f^\tau\right)_\tau=f^\tau_0$, thus $y$
(\ref{y_special}) can be recovered by:
\begin{equation}
\label{K1}
y_\tau=f^\tau_{0},\quad \tau=1,\ldots,T.
\end{equation}
And the potential can be found by
\begin{equation}
\label{K2}
b_n=\frac{y_{n+1}+y_{n-1}}{y_n}, n=1,\ldots,T-1.
\end{equation}

\subsection{Factorization method}

We make use the fact that matrix $C^T$ has a special structure -- it is a product of triangular matrix and its conjugate.
We rewrite the operator $W^T=\overline W^TJ$ as
\begin{equation*}
W^Tf=\begin{pmatrix}
1 & w_{1,1} & w_{1,2} & \ldots & w_{1,T-1}\\
0 & 1 & w_{2,2} &  \ldots & w_{2, T-1}\\
\cdot & \cdot & \cdot & \cdot & \cdot \\
0 & \ldots & 1& \ldots & w_{k, T-1}\\
\cdot & \cdot & \cdot  & \cdot & \cdot \\
0 & 0 & 0  & \ldots & 1
\end{pmatrix}
\begin{pmatrix}
0 & 0 & 0 & \ldots & 1\\
0 & 0 & 0 & \ldots  & 0\\
\cdot & \cdot & \cdot & \cdot &  \cdot \\
0 & \ldots & 1& 0  & 0\\
\cdot & \cdot & \cdot & \cdot &  \cdot \\
1 & 0 & 0 & 0 &  0
\end{pmatrix}
\begin{pmatrix} f_{0}\\
f_{2}\\
\cdot\\
f_{T-k-1}\\
\cdot\\
f_{T-1}
\end{pmatrix}
\end{equation*}
Using the definition (\ref{C_T_def}) and the invertibility of
$W^T$ (cf. Theorem \ref{teor_control}), we have:
\begin{equation*}
C^T=\left(W^T\right)^*W^T,\quad \text{or} \quad
\left(\left(W^T\right)^{-1}\right)^*C^T\left(W^T\right)^{-1}=I.
\end{equation*}
We can rewrite the latter equation as
\begin{equation}
\label{C_T_eqn_ker} \left(\left(\overline
W^T\right)^{-1}\right)^*\overline C^T\left(\overline
W^T\right)^{-1}=I,\quad \overline C^T=JC^TJ.
\end{equation}
Here the matrix $\overline C^T$ has the entries:
\begin{equation}
\label{C_overline_repr}
\overline C_{ij}=C_{T+1-j,T+1-i},\quad  \overline C^T=
\begin{pmatrix}
1 & r_1 & r_2 & \ldots & r_{T-1}\\
r_1 & 1+r_2 & r_1+r_3 & \ldots
&..\\
r_3 & r_1+r_3 & 1+r_2+r_4 & \ldots & ..\\
\cdot & \cdot & \cdot & \cdot & \cdot \\
\end{pmatrix},
\end{equation}
and operator $\left(\overline W^T\right)^{-1}$ has the form
\begin{equation}
\label{W_T_bar} \left(\overline W^T\right)^{-1}=\begin{pmatrix}
1 & \widetilde k_{11} & \widetilde k_{12}& \ldots & \widetilde k_{1,T-1} \\
0 & 1 & \widetilde k_{22} &\ldots &..\\
\cdot & \cdot & \cdot & \cdot & \widetilde k_{T-1,T-1} \\
0 &\ldots &\ldots & 0 & 1
\end{pmatrix},
\end{equation}
where $\widetilde k_{\alpha,\alpha}=-w_{\alpha,\alpha},$
$\alpha=1,\ldots,T-1.$ So we can rewrite (\ref{C_T_eqn_ker}) as
\begin{equation*}
\begin{pmatrix}
1 & 0 & . &  0 \\
k_{11} & 1 & 0  &.\\
\cdot & \cdot & \cdot & \cdot  \\
k_{T-1,1} & . & .  & 1
\end{pmatrix}
\begin{pmatrix}
\overline c_{11} & .. & .. &  \overline c_{1T} \\
.. & .. & ..  &..\\
\cdot & \cdot & \cdot & \cdot  \\
\overline c_{T1} &.. &   & \overline c_{TT}
\end{pmatrix}\begin{pmatrix}
1 & k_{11} & k_{21}&  .. \\
0 & 1 & k_{22}  &..\\
\cdot & \cdot & \cdot & \cdot  \\
0 &\ldots &\ldots  & 1
\end{pmatrix}=\begin{pmatrix}
1 & 0 & .. &  0 \\
0 & 1 & ..  &0\\
\cdot & \cdot & \cdot & \cdot  \\
0 &0 &  . & 1
\end{pmatrix}
\end{equation*}
In the above equation $\overline C_{ij}$ are given (see
(\ref{C_overline_repr})), the entries $k_{ij}$ of
$\left(\left(\overline W^T\right)^{-1}\right)^*$ are unknown. We
denote by $K_i:=\left(k_{i1},k_{i2},\ldots,
k_{ii},1,0,\ldots,0\right)$ the $(i+1)-$th row ($i=0,\dots,T-1$)
in $\left(\left(\overline W^T\right)^{-1}\right)^*$, then we have
\begin{equation*}
K_i\overline C^T K^*_j=\delta_{i,j}.
\end{equation*}
We use this equality in the form
\begin{equation}
\label{K_eqn} K_i\overline C^T K^*_j=0,\quad i<j.
\end{equation}
Notice that $K_0=\left(1,0,\ldots, 0\right)$. The second row $K_1$
can be recovered from $K_0\overline C^T K^*_1=0$, which is
equivalent to
\begin{equation}
\label{eq_1} \overline c_{11}k_{11}+\overline c_{21}=0,\,\,
\text{or}\,\, k_{11}=-\frac{\overline c_{21}}{\overline
c_{11}}=-\overline c_{21}.
\end{equation}
The third row $K_2$ we recover from the pair of equations
$K_0\overline C^T K^*_2=0,$ $K_1\overline C^T K^*_2=0$, which are
equivalent to
\begin{equation*}
\begin{pmatrix}
1 & 0\\
k_{11}& 1
\end{pmatrix}
\begin{pmatrix}
\overline c_{11} & \overline c_{12} & \overline c_{13}\\
\overline c_{21}& \overline c_{22} & \overline c_{23}
\end{pmatrix}
\begin{pmatrix}
k_{21} \\
k_{22}\\
1
\end{pmatrix}=\begin{pmatrix}
0\\
0
\end{pmatrix}.
\end{equation*}
Due to the invertibility of $\begin{pmatrix}
1 & 0\\
k_{1,1}& 1
\end{pmatrix}$ we can rewrite the latter equation as
\begin{equation}
\label{C_T_2sol}
\begin{pmatrix}
\overline c_{11} & \overline c_{12}\\
\overline c_{21}& \overline c_{22}
\end{pmatrix}
\begin{pmatrix}
k_{21} \\
k_{22}\\
\end{pmatrix}=-\begin{pmatrix}
\overline c_{13}\\
\overline c_{23}
\end{pmatrix}.
\end{equation}
We introduce the notation, by $c_i^k$ we denote the i-th column in
the matrix $\overline C^T$ truncated by first $k$ elements:
\begin{equation*}
\overline c_i^k:=\begin{pmatrix} \overline c_{1i} & \overline
c_{2i}&\ldots & \overline c_{ki}
\end{pmatrix}^*.
\end{equation*}
Since $C^T$ is invertible, (\ref{C_T_2sol}) has a unique solution,
moreover, we can infer that
\begin{equation*}
\label{eq2} k_{22}=-\frac{\det\begin{pmatrix}
\overline c_{11} & \overline c_{13}\\
\overline c_{21}& \overline c_{23}
\end{pmatrix}}{\det\begin{pmatrix}
\overline c_{11} & \overline c_{12}\\
\overline c_{21}& \overline c_{22}
\end{pmatrix}}=-\det(\overline c_{1}^2,\overline c_{3}^2).
\end{equation*}
Assume that we have already recovered $K_0,K_1,\ldots,K_l$, to
recover $K_{l+1}$ we need to consider the equations $K_0\overline
C^T K^*_{l+1}=0,$ $K_1\overline C^T K^*_{l+1}=0,\ldots,$
$K_l\overline C^T K^*_{l+1}=0$, which are equivalent to
\begin{equation*}
\begin{pmatrix}
1 & 0 & ..& 0\\
k_{11}& 1& 0&..\\
.&.&.&.\\
k_{l1}& k_{l2}&.&1
\end{pmatrix}
\begin{pmatrix}
\overline c_{11} & ..& .. & \overline c_{1,l+2}\\
.. &.. &.. &.. \\
.. &.. &.. &.. \\
\overline c_{l+1,1}& .. & ..& \overline c_{l+1,l+2}
\end{pmatrix}
\begin{pmatrix}
k_{l+1,1} \\
k_{l+1,2}\\
..\\
1
\end{pmatrix}=\begin{pmatrix}
0\\
..\\
..\\
0
\end{pmatrix}.
\end{equation*}
We can rewrite the latter equation as
\begin{equation}
\label{GL_eqn}
\begin{pmatrix}
\overline c_{1,1} & ..& .. & \overline c_{1,l+1}\\
.. &.. &.. &.. \\
.. &.. &.. &.. \\
\overline c_{l+1,1}& .. & ..& \overline c_{l+1,l+1}
\end{pmatrix}\begin{pmatrix}
k_{l+1,1} \\
k_{l+1,2}\\
..\\
k_{l+1,l+1}
\end{pmatrix}+\begin{pmatrix}
\overline c_{1,l+2} \\
\overline c_{2,l+2}\\
..\\
\overline c_{l+1,l+2}
\end{pmatrix}=0
\end{equation}
Due to the invertibility of $C^T$ the latter equation has unique
solution, moreover
\begin{equation}\label{k_rec}
k_{l+1,l+1}=-\det(\overline c_{1}^{l+1},\overline
c_{2}^{l+1},\ldots \overline c_{l}^{l+1},\overline
c_{l+2}^{l+1}),\quad l=0,\ldots, T-2.
\end{equation}
Having recovered $k_{\alpha,\alpha}=-w_{\alpha,\alpha}$, we
recover the potential by (see (\ref{Goursat}))
\begin{equation}
\label{Pot_rec} b_n=w_{n-1,n-1}-w_{n,n}=k_{n,n}-k_{n-1,n-1},\quad
n=1,\ldots,T-1.
\end{equation}

\subsection{Gelfand-Levitan equations}

If we introduce $\widetilde C^T$ by
\begin{equation}
\label{widetild_def}
\overline C^T=I+\widetilde C^T,
\end{equation}
(see (\ref{C_T_repr}),(\ref{C_overline_repr})), then we can
rewrite (\ref{GL_eqn}) for $l=T-2$ as
\begin{equation*}
\left(I+\widetilde C^T\right)K_T+\widetilde C_{T}=0,
\,\,\text{where}\,\, K_T=\begin{pmatrix} k_{T-1,1}\\
k_{T-1,2}\\
.\\
k_{T-1,T-1}
\end{pmatrix},\,\,\widetilde C_{T}=\begin{pmatrix} \widetilde C^T_{1,T}\\
\widetilde C^T_{2,T}\\
.\\
\widetilde C^T_{T-1,T}
\end{pmatrix}
\end{equation*}
or as a system
\begin{equation}
\label{GL_eqn_fin} k_{T-1,\beta}+\sum_{j=1}^{T-1}\widetilde
C^T_{\beta,j} k_{T-1,j}+\widetilde C^T_{\beta,T}=0,\quad
\beta=1,\ldots,T-1.
\end{equation}
If we  pass to (more standard) entries of $\left(\overline
W^T\right)^{-1}$
\begin{equation}\label{K_tilde}
\widetilde k_{\alpha,\beta}=k_{\beta,\alpha},
\end{equation}
then ($\ref{GL_eqn_fin}$) can be rewritten as
\begin{equation}
\label{GL_eqn_fin1} \widetilde
k_{\beta,T-1}+\sum_{j=1}^{T-1}\widetilde C^T_{\beta,j} \widetilde
k_{j,T-1}+\widetilde C^T_{\beta,T}=0,\quad \beta=1,\ldots,T-1.
\end{equation}
The last equation is an analogue of Gelfand-Levatan equation for
continuous problem \cite{AM,BM_1}. We conclude this section with
\begin{theorem}
The kernel of the operator $\left(\overline W^T\right)^{-1}$ (see
(\ref{W_T_bar})) satisfies equation (\ref{GL_eqn_fin1}), where the
entries $\widetilde C^T_{i,j}$ are defined in
(\ref{widetild_def}), (\ref{C_T_repr}).
\end{theorem}
The equation in (\ref{GL_eqn_fin}) has a unique solution due to
the invertibility of $C^T$. The potential can be recovered by
(\ref{Pot_rec}).

Now we make some remarks on the dependence of the connecting
operator $C^T$ and the solution of the inverse problem equations
(i.e. the potential) on the inverse data. As a direct consequence
of (\ref{C_T_repr}) we can formulate the folowing
\begin{remark}
The operator $C^T$ depends on $R^{2T-2}$, i.e. it depends on the
potential $(b_1,\ldots,b_{T-1})$, so the results obtained from
$C^T$ via Krein-type equations (\ref{C_T_Krein}),
(\ref{K1}),(\ref{K2}), factorization method (\ref{k_rec}),
(\ref{Pot_rec}) and Gelfand-Levitan type equations
(\ref{GL_eqn_fin1}), (\ref{Pot_rec}) are best possible.
\end{remark}

In the subsection on the factorization method we used the fact
that $\det C^\tau=1,$ $\tau=2,\ldots, T$. More precisely, we used
it in the form $\det(\overline c_{1}^{\tau},\overline
c_{2}^{\tau},\ldots \overline c_{\tau}^{\tau})=1$. That fact
actually says that not all elements in the response vector are
independent. Indeed: the element $k_{11}$ we recovered (see
(\ref{eq_1})) from $\overline c_{21}$, i.e. from $r_1$. The
element $k_{22}$ we recovered from $\overline c_{11}$, $\overline
c_{13}$, $\overline c_{21}$, $\overline c_{23}$, that is from
$r_1,$ $r_2,$ $r_3$. But since $\det(\overline c_{1}^{2},\overline
c_{2}^{2})=1$, we have that $r_2=r_1^2$, so in fact $k_{22}$ was
recovered from $r_1$ and $r_3.$ Arguing in the same fashion we see
that $r_{2k}$ depends on $r_{2l+1}$, $l=0,\ldots,k-1$. So we
recovered $(k_{11},\ldots, k_{T-1,T-1})$ from the response vector
$(r_0,r_1,\ldots,r_{2t-2})$, $r_0=1$, whose components with even
numbers depend in explicit form on the components with odd
numbers. That observation plays an important role in the next
subsection.

\subsection{Characterization of the inverse data.}
In the second section we considered the forward problem
(\ref{Jacobi_dyn}), for the potential $(b_1,\dots, b_{T-1})$ we
constructed the matrix $W^T$ (\ref{Jac_sol_rep}), (\ref{Goursat}),
the response vector $(1,r_1,\dots,r_{2T-2})$ (see (\ref{R_def}))
and the connecting operator $C^T$ by formula (\ref{C_T_repr}). It
will be more convenient for us to deal with the rotated matrix
$\overline C^T$ defined in (\ref{C_overline_repr}). From the
representation $\overline C^T= (\overline W^T)^* \overline W^T$
and triangularity of $\overline W^T$ we know that
\begin{equation*}
\det\overline C^l=1\quad\forall l=1,\dots,T.
\end{equation*}
Also we have proved that if coefficients $r_1,\dots,r_{2T-2}$
correspond to some potential $(b_1,\dots,b_{T-1})$ then we can
recover the potential using (\ref{k_rec})-(\ref{Pot_rec}).

Now we set up a question: can one determine whether a vector
$(1,r_1,r_2,\ldots,r_{2T-2})$ is a response vector for the
dynamical system (\ref{Jacobi_dyn}) with a potential
$(b_1,\ldots,b_{T-1})$ or not? The answer is the following
theorem.
\begin{theorem}
The  vector $(1,r_1,r_2,\ldots,r_{2T-2})$ is a response vector for
the dynamical system (\ref{Jacobi_dyn}) if and only if the matrix
$C^T$ (\ref{C_T_repr}) is positively definite and $\det C^l=1,$
$l=1,\ldots,T$.
\end{theorem}
\begin{proof}
First we observe that in the conditions of the theorem we can
substitute $C^T$ by $\overline C^T$ (\ref{C_overline_repr}). The
necessary part of the theorem is proved in the preceding sections.
We are left to prove the sufficiency of these conditions.

Let we have a vector $(1,r_1,\dots,r_{2T-2})$ such that the matrix
$\overline C^{T}$ constructed from it using
(\ref{C_overline_repr}) satisfies conditions of the theorem. Then
we can construct the potential $(b_1,\dots,b_{T-1})$ using
(\ref{k_rec})-(\ref{Pot_rec}) and consider the dynamical system
(\ref{Jacobi_dyn}) with this potential. For this system we
construct the connecting operator $C_{new}^{T}$ and its rotated
$\overline C^{T}_{new}$ using (\ref{Goursat}), (\ref{R_def}),
(\ref{C_T_repr}) and (\ref{C_overline_repr}). We will show that
the matrices $\overline C^{T}$ and $\overline C^{T}_{new}$
coincide.

First of all we note that we have two matrices constructed by
(\ref{C_overline_repr}), one comes from the vector
$(1,r_1,\dots,r_{2T-2})$ and the other comes from
$(1,r^{new}_1,\dots,r^{new}_{2T-2})$. Also they have a common
property that $\det \overline C^l=\det \overline C^l_{new}=1$ for
all $l=1,\dots,T$ (one by theorem's condition and the other by
representation $\overline C^T_{new}= (\overline W^T_{new})^*
\overline W^T_{new}$).

Secondly we note that if we calculate the potential
$(b_1,\dots,b_{T-1})$ using (\ref{k_rec})-(\ref{Pot_rec}) from any
of $\overline C^{T}$ and $\overline C^{T}_{new}$ matrices, we get
the same answer.

Therefore we have two matrices of the type (\ref{C_overline_repr})
with the unit principal minors and the property
\begin{equation}\label{c_eq_c}
\det(\overline c_{1}^{l+1},\ldots \overline c_{l}^{l+1},\overline
c_{l+2}^{l+1})  =\det(\overline {c_{new}}_{1}^{l+1},\ldots
\overline {c_{new}}_{l}^{l+1},\overline {c_{new}}_{l+2}^{l+1})
\quad \forall l=0,\dots,T-2.
\end{equation}
If we  look at (\ref{c_eq_c}) for $l=0$, we see that $r_1=
r_1^{new}$. From the fact the for both matrices $\overline C^T,$
$\overline C^T_{new}$ the principal minors of the second order are
equal to one, we infer that $r_2=r_2^{new}$. We continue this
procedure, and from (\ref{c_eq_c}) with $l=n$ we infer that
$r_{2n+1}=r_{2n+1}^{new}$ and from equality to one of principal
minor of the order $n+2$ of $\overline C^T,$ $\overline
C^T_{new}$, we can infer that $r_{2n+2}=r_{2n+2}^{new}$ for all
$n=2,\dots,T-2$ by induction. That finishes the proof.
\end{proof}

\section{Spectral representation of $C^T$ and $r_t$.}

In this section we consider the inverse spectral problem and show
the connection of the spectral (\ref{Ortog}), (\ref{SP_data}) and
dynamical (\ref{R_def}), (\ref{R_def_int}) inverse data. If we
introduce the special control $\delta=(1,0,0,\ldots)$, then the
kernel of response operator (\ref{R_def_int}) is
\begin{equation}
\label{con1} r_t^i=\left(R_i\delta\right)_t=v^\delta_{1,t},
\end{equation}
on the other hand, we can use (\ref{Jac_sol_rep_int}),
(\ref{c_koeff_repr}) to obtain:
\begin{equation}
\label{con2}
v^\delta_{1,t}=\sum_{k=1}^N\frac{1}{\rho_k}T_t(\lambda_k).
\end{equation}
So on introducing the spectral function
\begin{equation}
\label{Spectr_fun}
\rho^N(\lambda)=\sum_{\{k\,|\,\lambda_k<\lambda\}}\frac{1}{\rho_k},
\end{equation}
from (\ref{con1}), (\ref{con2}) we deduce that
\begin{equation*}
r_t^i=\int_{-\infty}^\infty T_t(\lambda)\,d\rho^N(\lambda),\quad
t\in \mathbb{N}.
\end{equation*}
Let us evaluate $(C^T_if,g)$ for $f,g\in \mathcal{F}^T$, using the
expansion (\ref{Jac_sol_rep_int}):
\begin{eqnarray*}
(C^T_if,g)=\sum_{n=1}^N
v^f_{n,T}v^g_{n,T}=\sum_{n=1}^N\sum_{k=1}^N
\frac{1}{\rho_k}T_T\left(\lambda_k\right)*f\varphi^k_n \,
\sum_{l=1}^N
\frac{1}{\rho_l}T_T\left(\lambda_l\right)*g\varphi^l_n\\
=\sum_{k=1}^N\frac{1}{\rho_k}T_T(\lambda_k)*f
T_T(\lambda_k)*g=\int_{-\infty}^\infty
\sum_{l=0}^{T-1}T_{T-l}(\lambda)f_l
\sum_{m=0}^{T-1}T_{T-m}(\lambda)g_m\,d\rho^N(\lambda)
\end{eqnarray*}
from the equality above it is evident that (cf. (\ref{C_T_repr}))
\begin{equation}
\{C^T_i\}_{l+1,m+1}=\int_{-\infty}^\infty
T_{T-l}(\lambda)T_{T-m}(\lambda)\,d\rho^N(\lambda), \quad
l,m=0,\ldots,T-1.
\end{equation}
Let us consider the spectral problem
\begin{equation}
\label{Spec_sol1} \left\{
\begin{array}l \phi_{i+1}+\phi_{i-1}-b_n\phi_i=\lambda\phi_i,\quad n=0,\ldots,N+1,\\
\phi_0=0,\,\,\phi_{N+1}=0.
\end{array}
\right.
\end{equation}
In the second section we construct the spectral data for this
problem -- eigenvalues of the corresponding Hamiltonian and
norming coefficients (\ref{Ortog}), (\ref{SP_data}). Now we answer
the question how to recover the potential $(b_1,\ldots,b_N)$ from
this data.

Our strategy will be to use the dynamical approach from the forth
section to treat this IP. First we observe that to know
(\ref{SP_data}) is the same as to know the spectral function
(\ref{Spectr_fun}). Consider the system (\ref{Jacobi_dyn}) with
the same potential $b_n$ for $n=1,\ldots,N$. We notice that as
explained in the beginning of the section four, $R^{2N}=R^{2N}_i$
and correspondingly, $r_t=r_t^i,$ $t=1,\ldots,2N.$ Due to this, we
deduce that $C^T=C^T_i$ for $T=N+1$. Thus, the inverse problem can
be solved in the following way: from the spectral data
(\ref{SP_data}) we construct the spectral function by
(\ref{Spectr_fun}). Then we construct
\begin{eqnarray*}
r_t=r^i_t=\int_{-\infty}^\infty
T_t(\lambda)\,d\rho^N(\lambda),\quad t=1,\ldots,2N,\\
C^T_{lm}=\{C^T_i\}_{l+1,m+1}=\int_{-\infty}^\infty
T_{T-l}(\lambda)T_{T-m}(\lambda)\,d\rho^N(\lambda), \quad
l,m=0,\ldots,N-1.
\end{eqnarray*}
After we have in hands the connecting operator, we can use the
methods of section four to find $(b_1,\ldots,b_N).$

\noindent{\bf Acknowledgments}

The research of Victor Mikhaylov was supported in part by NIR
SPbGU 11.38.263.2014 and RFBR 14-01-00535. Alexandr Mikhaylov was
supported by RFBR 14-01-00306; A. S. Mikhaylov and V. S. Mikhaylov
were partly supported by VW Foundation program "Modeling,
Analysis, and Approximation Theory toward application in
tomography and inverse problems."


\begin{thebibliography}{99}

\bibitem{AMM}
{S. A. Avdonin, A. S. Mikhaylov, V. S. Mikhaylov,} \textit{On some
applications of the Boundary Control method to spectral estimation
and inverse problems,} Nanosystems: Phys. Chem. Math. {\bf 6} (1),
(2015). 63--78.

\bibitem{AM}
{S. A. Avdonin, V. S. Mikhaylov,} \textit{The boundary control
approach to inverse spectral theory,} Inverse Problems {\bf 26},
2010, no. 4, 045009, 19 pp.

\bibitem{B07}
{M. I. Belishev}, \textit{Recent progress in the boundary control
method}, Inverse Problems, {\bf 23} (5), (2007). R1.
doi:10.1088/0266-5611/23/5/R01

\bibitem{B13}
{M. I. Belishev}, \textit{C*-Algebras in reconstruction of
manifolds}, Nanosystems: Phys. Chem. Math., {\bf 4} (4), (2013),
484-489.

\bibitem{BM_1}
{M.I.Belishev, V.S.Mikhailov}. \textit{Unified approach to
classical equations of inverse problem theory.} {Journal of
Inverse and Ill-Posed Problems}, {\bf 20} no. 4, (2012), 461--488.

\end{thebibliography}
\end{document}